\documentclass[10pt]{amsart}
\usepackage{latexsym}
\input diagrams.tex

\def\sw#1{{\sb{(#1)}}}
\def\su#1{{\sp{(#1)}}} 
\def\tens{\mathop{\otimes}}

\def\<{{\langle}}
\def\>{{\rangle}}

\def\eps{\epsilon}

\def\note#1{{}}

\def\note#1{}
\def\M{{\bf M}}
\def\N{{\mathbb N}}

\def\cM{{\mathfrak M}}

\def\cN{{\mathfrak N}}
\def\cC{{\mathfrak C}}

\def\cD{{\mathfrak D}}

\def\lhom#1#2#3{{{}\sb{#1}{\rm Hom}(#2,#3)}}
\def\rhom#1#2#3{{{\rm Hom}\sb{#1}(#2,#3)}}

\def\Rrhom#1#2#3#4{{{\rm Hom}\sp{#1}\sb{#2}(#3,#4)}}

\def\beq{\begin{equation}}
\def\eeq{\end{equation}}
\def\DC{{\Delta_\cC}}
\def \eC{{\eps_\cC}}
\def\DD{{\Delta_\cD}}
\def \eD{{\eps_\cD}}

\def\ff{{\mathfrak f}}

\def\con{{\bf Conn}}

\def\bdi{\begin{diagram}}
\def\edi{\end{diagram}}

\def\Label#1{\label{#1}\ifmmode\llap{[#1] }\else 
\marginpar{\smash{\hbox{\tiny [#1]}}}\fi}
\def\Label{\label}

\newtheorem{proposition}{Proposition}[section]

\newtheorem{theorem}[proposition]{Theorem}

\theoremstyle{definition}

\newtheorem{example}[proposition]{Example}

\theoremstyle{remark}
\newtheorem{remark}[proposition]{Remark}

\newcounter{c}

\newcommand{\etyk}[1]{\vspace{-7.4mm}$$\begin{equation}\Label{#1}
\addtocounter{c}{1}}
\renewcommand{\]}{\ifnum \value{c}=1 $$\else \end{equation}\fi}
\begin{document}

\title{The structure of corings with a grouplike element}
\author{Tomasz Brzezi\'nski}
\address{Department of Mathematics, University of Wales Swansea,
Singleton Park, Swansea SA2 8PP, U.K.}
\email{T.Brzezinski@swansea.ac.uk}
\urladdr{http//www-maths.swan.ac.uk/staff/tb}
\begin{abstract}
Characteristic properties of corings with a grouplike 
element are analysed. Associated differential
graded rings are studied. A correspondence between 
categories of comodules and flat
connections is established. A generalisation of the Cuntz-Quillen 
theorem relating existence of connections in a module to projectivity 
of this module is proven.
\end{abstract}
\maketitle

\section{Introduction}
The notion of a coring was introduced by Sweedler in \cite{Swe:pre} as a
generalisation of coalgebras over commutative rings to the case of
non-commutative rings. While the theory of coalgebras has been developed
quite substantially over the last quarter of the century, most notably
after the discovery of large classes of examples coming from quantum 
groups, the progress in the coring theory was hampered by the lack of
examples. Recently M.\ Takeuchi realised that the 
compatibility condition between an algebra and a coalgebra known as an 
entwining \cite{BrzMa:coa} can be recast in terms of a coring. 
Entwining structures and the associated categories of entwined modules
\cite{Brz:mod} unify and generalise various  categories of Hopf modules
such as relative Hopf modules \cite{Tak:cor}, \cite{Doi:str},
Yetter-Drinfeld modules \cite{Yet:rep}, \cite{RadTow:yet} or
Doi-Koppinen Hopf modules \cite{Doi:uni}, \cite{Kop:var}. Several
classes of examples of entwining structures are known (cf.\ recent
monograph \cite{CaeMil:gen}). Thus it turns out
that the theory of corings is very rich in examples. Once a coring
is built from an entwining structure, the category of entwined modules
is isomorphic to the category of comodules of this coring. Various
properties of entwined modules can be understood and more simply 
derived from general properties of corings and their comodules. A
programme of studying of corings from this point of view was initiated
in \cite{Brz:str} and is now  gaining a momentum with several new
results reported in the area (cf.\ \cite{Abu:rat}, \cite{Brz:tow},
\cite{ElKGom:sem}, \cite{Gom:sep}, \cite{Wis:wea},
\cite{Wis:com}). 

The present paper is devoted to the study of general properties of those
corings, which have a grouplike
element. These are corings most closely related to the original example
studied by Sweedler in \cite{Swe:pre}, and they reveal the richness of 
structure, with vistas on the (noncommutative) descent theory and
 noncommutative geometry. We begin
in Section~2 by recalling the basic definitions and examples, and by
deriving elementary properties of corings with a grouplike element. We
then proceed in Section~3 to 
recall their connection with graded differential rings, 
and reveal
their close relationship with noncommutative connections in Section~4.

All rings in this paper are associative and unital, and the unit of a ring
$R$ is denoted by $1_R$, while $\M_R$  denotes the category of right
$R$-modules, and ${}_R\M$ the left $R$-modules etc. All algebras and
coalgebras are over a commutative ring $k$, undecorated tensor products
and homomorphisms are over $k$. The identity
morphism for
any object $X$ in any category is denoted by the same symbol $X$.

\section{Basic properties}
\Label{sec.basic.grouplike}
Given a ring $R$, an $(R,R)$-bimodule $\cC$ is called an $R$-coring if
there exist $(R,R)$-bimodule maps $\Delta_\cC:\cC\to
\cC\otimes_R\cC$  and $\eps_\cC:\cC\to R$, such that  
$$ 
(\Delta_\cC\otimes_R\cC)\circ\Delta_\cC = 
    (\cC\otimes_R\Delta_\cC)\circ \Delta_\cC, \quad 
(\eps_\cC\otimes_R\cC)\circ\Delta_\cC = 
(\cC\otimes_R\eps_\cC)\circ \Delta_\cC = 
\cC.
$$ 
The map $\Delta_{\cC}$ is called the {\em coproduct}\index{coproduct} 
while the map 
$\eps_{\cC}$ is called the {\em counit}\index{counit} of $\cC$. A
morphism of $R$-corings $\cC$ and $\cD$ is an $(R,R)$-bimodule map $\ff:
\cC\to \cD$ such that $\DD\circ\ff = (\ff\otimes_R\ff)\circ\DC$ and
$\eD\circ\ff = \eC$. Given an
$R$-coring $\cC$, a right $\cC$-comodule is a right $R$-module $\cM$
together with 
a right $R$-module map $\rho^\cM :\cM\to \cM\otimes_R\cC$ such that
$$
(\rho^\cM\otimes_R\cC)\circ\rho^\cM = 
    (\cM\otimes_R\Delta_\cC)\circ \rho^\cM, \quad 
(\cM\otimes_R\eps_\cC)\circ \rho^\cM = 
\cM.
$$ 
A map between  right $\cC$-comodules is a right $R$-module map
respecting the coactions, i.e., $\ff:\cM\to \cN$ satisfies
$\rho^\cN\circ \ff = (\ff\otimes_R\cC)\circ\rho^\cM$. The category of
right $\cC$-comodules is denoted by $\M_R^\cC$, and their morphisms by
${\rm Hom}_R^\cC(-,-)$. The actions of $\DC$ and $\rho^\cM$ on 
elements are denoted by the Sweedler notation $\DC(c) = c\sw 1\otimes_R
c\sw 2$, $\rho^\cM(m) = m\sw 0\otimes_R m\sw 1$. Similarly one defines
left $\cC$-comodules, and bicomodules.

Any ring $R$ can be viewed as an $R$-coring via the identity maps. In
this case the category of comodules is the same as the category of right
$R$-modules. Furthermore, given a ring extension $S\to R$, the
$(R,R)$-bimodule $\cC = R\otimes_S R$ is an $R$-coring with the
coproduct $\DC(r\otimes_S r') = r\otimes_S 1_R\otimes_S r'$ and the
counit $\eC(r\otimes_s r') = rr'$, for all $r,r'\in R$. We refer to this
coring as the {\em canonical Sweedler's coring} associated to a ring
extension $S\to R$. The category of $\cC$-comodules in this case is
isomorphic to the category of (noncommutative) descent data \cite{KnuOja:des},
\cite{Cip:dis} for the ring extension $S\to R$. Finally, let 
$(A,C)_\psi$ be an entwining structure
over a commutative ring $k$, i.e., $A$ is a $k$-algebra, $C$ is a
$k$-coalgebra with a counit $\eps$ and coproduct $\Delta$, and $\psi:
C\otimes A\to A\otimes C$ is a $k$-linear map satisfying mixed
distributive law conditions (cf.\ \cite{BrzMa:coa} or \cite{Brz:str} for
details). Then $\cC = A\otimes C$ is an $A$-coring with the following
structure. Left $A$-multiplication is given by the product in $A$, 
right
$A$-multiplication is provided by the map $\psi$, i.e., $(a\otimes c)a'
= a\psi(c\otimes a')$, the coproduct and the counit are $\DC =
A\otimes \Delta$ and $\eC= A\otimes\eps$. In this case $\M^\cC_R =
\M^A_C(\psi)$, the latter denoting the category of $(A,C)_\psi$-entwined
modules, i.e., right $A$-modules and $C$-comodules $M$ with a compatibility
condition $\rho^M(ma) = m\sw 0\psi(m\sw 1\otimes a)$, for all $m\in M,
a\in A$, introduced in \cite{Brz:mod}. 

A non-zero element $g$ of an $R$-coring $\cC$ 
is said to be a {\em grouplike} element
 if $\Delta_{\cC}(g) = 
g\otimes_{R}g$ and  
$\eps_{\cC}(g) = 1_{R}$. A ring $R$ viewed as a trivial $R$-coring has a
grouplike element $1_R$. Similarly, the canonical Sweedler's coring
associated to $S\to R$ has a grouplike element $1_R\otimes_S1_R$. In
general the existence of a grouplike element determines when a ring is
itself a comodule of a coring (cf.\ \cite[Lemma~5.1]{Brz:str}).
Precisely, an 
$R$-coring $\cC$ has a grouplike element $g$ if and only if $R$ is 
a right  $\cC$-comodule. The coaction reads $\rho^R(r) = gr$. For
example, in the case of an entwining structure $(A,C)_\psi$, the
corresponding coring has a grouplike element if and only if $A$ is an
entwined module, In this case $g = \rho^A(1_A)\in A\otimes C$, where
$\rho^A$ is a coaction of $C$ on $A$.

Given an $R$-coring with a grouplike element $g$ and a right
$\cC$-comodule $\cM$ with a coaction $\rho^\cM: \cM\to\cM\otimes_R 
\cC$, one defines {\em
$g$-coinvariants} 
of $\cM$ as an Abelian group
$$
\cM^{co\cC}_g = \{m\in\cM \;|\; \rho^\cM(m) = m\otimes_R g\}.
$$
The $g$-coinvariants $S$ of $R$ form a subring
 of $R$ equal to the
centraliser of $g$ in $R$, i.e.,
$S = \{s\in R, \; | \; sg=gs\}$.

Given any right $S$-module $M$, where $S$ are $g$-coinvariants of $R$,
one can consider $M\otimes_SR$ as a right $\cC$-comodule via the
coaction
$$
\rho^{M\otimes_SR}: M\otimes_SR\to M\otimes_SR\otimes_R\cC \cong
M\otimes_S\cC, \qquad m\otimes_Sr\mapsto m\otimes_Sgr.
$$
The assignment $M\mapsto M\otimes_S R$ defines a covariant functor
$-\otimes_SR:\M_S\to \M^\cC_R$, known as an {\em induction 
functor}\index{induction functor}. If $f:M\to N$ is a morphism in 
$\M_{S}$ then $f\otimes_{S}R: m\otimes_{S}r\mapsto f(m)\otimes_{S}r$. 
Note that $f\otimes_{S}R$ is a morphism in $\M^{\cC}_{R}$ since 
$\rho^{N\otimes_{S}R}(f(m)\otimes_{S}r) = f(m)\otimes_{S}gr$ and 
$(f\otimes_{S}R\otimes_{R}\cC)\circ 
\rho^{M\otimes_{S}R}(m\otimes_{S}r) = 
(f\otimes_{R}\cC)(m\otimes_{R}gr) = f(m)\otimes_{S}gr$.
In the opposite direction \cite[Proposition~5.2]{Brz:str},  
the assignment 
$G_{g}:\M^\cC_R\to \M_S$, $\cM\mapsto \cM^{co\cC}_g$,
defines a covariant functor, known as a {\em $g$-coinvariants functor},
which is the right adjoint of the 
induction
functor $-\otimes_SR:\M_S\to \M^\cC_R$. On morphisms $G_{g}$ acts as a 
restriction of the domain, i.e., for any $\ff:\cM\to \cN$ in 
$\M_{R}^{\cC}$, $G_{g}(\ff) = \ff\mid_{\cM_{g}^{co\cC}}$.

\begin{example}
Let $\cC$ be an $R$-coring with a grouplike element $g$, 
and let $S= R^{co\cC}_g$ be the
subring of $g$-coinvariants. View $\cC$ as a right or left
$\cC$-comodule via the coproduct $\DC$. Then $R\cong \cC^{co\cC}_g$ as
$(R,S)$-bimodules.
\Label{ex.coinv.C}
\end{example}
\begin{proof}
First note that since $\DC$ is an $(R,R)$-bimodule map, the set of
$g$-coinvariants $\cC^{co\cC}_g =\{c\in\cC\; |\; \DC(c) = c\otimes_R
g\}$ is a left $R$-module. Thus $\cC^{co\cC}_g$ is an $(R,S)$-bimodule.
Consider a left $R$-module map $\phi: R\to \cC_g^{co\cC}$, 
$r\mapsto rg$.
Since $S$ centralises $g$ in $R$, we have for all $s\in S$ and $r\in R$,
$\phi(rs) = rsg = rgs = \phi(r)s$,
i.e., $\phi$ is an $(R,S)$-bimodule map. It is well-defined since
$\DC(rg) =rg\otimes_Rg$, i.e., $rg\in \cC_g^{co\cC}$ for all $r\in R$,
as needed. Now consider
$\phi^{-1}= \eC\mid_{\cC^{co\cC}_g} :\cC^{co\cC}_g\to R$.
Clearly $\phi^{-1}$ is an $(R,S)$-bimodule map (since $\eC$ is an
$(R,S)$-bimodule map). Furthermore for all $r\in R$,
$\phi^{-1}\circ\phi(r) = \eC(rg) =r$. Notice that if $c\in
\cC^{co\cC}_g$, then $\DC(c) = c\otimes_R g$, hence applying
$\eC\otimes_R \cC$ we have $c=\eC(c)g$. Using this fact take any
$c\in \cC^{co\cC}_g$ and compute, $
\phi\circ\phi^{-1}(c) = \phi(\eC(c)) =\eC(c) g = c$.
This proves that $\phi^{-1}$ is the inverse of $\phi$ in ${}_R\M_S$,
i.e., $R\cong \cC^{co\cC}_g$ as $(R,S)$-bimodules. 
\end{proof}

The
$g$-coinvariants functor can be viewed as a hom-functor.
\begin{proposition}
Let $\cC$ be an $R$-coring with a grouplike element $g$ and let 
 $S= R^{co\cC}_g$ be the
subring of $g$-coinvariants. 
For any $\cM\in\M^\cC_R$ view $\Rrhom \cC RR\cM$ as a right
$S$-module via $(\ff s)(r) = \ff(sr)$, for all $\ff\in \Rrhom \cC
RR\cM$, $r\in R$ and $s\in S$. Then $\Rrhom \cC RR\cM\cong
\cM^{co\cC}_g$ as right $S$-modules.
\Label{prop.coinv.hom}
\end{proposition}
\begin{proof}
First note that  $\Rrhom \cC RR\cM$  is a right $S$-module as
stated. Indeed, for all $s\in S$, $r\in R$ and $\ff\in \Rrhom \cC RR\cM$
one computes
\begin{eqnarray*}
(\ff s)(r\sw 0)\otimes_R r\sw 1 &=& (\ff s)(1_R)\otimes_R gr =
\ff(s)\otimes_R gr = \ff(1_R)\otimes_R sgr\\
&=& \ff(1_R)\otimes_R gsr = \ff((sr)\sw 0)\otimes_R (sr)\sw 1 \\
& = &
\rho^\cM(\ff(sr)) = \rho^\cM((\ff s)(r)),
\end{eqnarray*}
i.e., $\ff s$ is a right $\cC$-comodule map as required. Clearly the
multiplication $\ff s$ is associative.
Now to any $\ff\in\Rrhom \cC RR\cM$ assign $m_\ff = \ff(1_R)\in \cM$.
Since $\ff$ is a right $\cC$-comodule map we have
$$
\rho^\cM(m_\ff) = \rho^\cM(\ff(1_R)) = \ff(1\sw 0)\otimes_R 1\sw 1 =
\ff(1_R)\otimes_R g = m_\ff\otimes_R g,
$$
so that $m_\ff \in \cM^{co\cC}_g$. Furthermore for all $s\in S$ we have
$ m_{\ff s} = (\ff s)(1_R) = \ff(s) = \ff(1_R) s = m_\ff s$. All 
this means that
the assignment $\ff\mapsto m_\ff$ is a right $S$-module morphism $\Rrhom \cC
RR\cM\to \cM^{co\cC}_g$.

Conversely to any $m\in \cM^{co\cC}_g$ assign $\ff_m\in \rhom RR\cM$,
$\ff_m : r\mapsto mr$. Since $m$ is in $g$-coinvariants, and $\rho^\cM$
is a right $R$-module morphism one can take
any $r\in R$ and compute
$$
\rho^\cM(\ff_m(r)) = \rho^\cM(mr) = m\otimes_R gr = \ff_m(1_R)\otimes
_R gr =\ff_m(r\sw 0)\otimes_R r\sw 1,
$$
i.e., $\ff_m$ is a right $\cC$-comodule map. Since the assignments
$\ff\mapsto m_\ff$ and $m\mapsto \ff_{m}$ 
described above are restrictions of the natural isomorphism $\rhom
RR\cM\cong \cM$, they are inverses to each other as required.
\end{proof}

The existence of a grouplike element $g$ in $\cC$ allows for a direct 
sum decomposition of $\cC$, 
 $\cC \cong R\oplus \ker\eps_{\cC}$ as $(S,R)$-bimodules  and 
 $(R,S)$-bimodules, where $S=R^{co\cC}_{g}$. The $(S,R)$-bimodule  
 isomorphism is given by $u_{R}:  c\mapsto 
 (\eps_{\cC}(c),g\eps_{\cC}(c)-c)$, with the inverse 
 $u_{R}^{-1}(r,c) = gr-c$, while the  $(R,S)$-bimodule isomorphism is 
 $u_{L}:  c\mapsto    (\eps_{\cC}(c),c- \eps_{\cC}(c)g)$, with the 
 inverse $u_{L}^{-1}(r,c) = c+rg$.
Since $\ker\eC$ is an $(R,R)$-bimodule there is a natural ring 
structure on $R\oplus \ker\eps_{\cC}$ 
with the product given by $(r,c)(r',c') = (rr', rc'+cr')$, 
for all $r,r'\in R$ and $c,c'\in\ker\eC$,
 and the unit $(1,0)$. Using the isomorphisms $u_{L}$ or $u_{R}$ this 
 ring structure can be pulled back to $\cC$, and the product comes out  as 
    $$
    cc' = \eps_{\cC}(c)c'+c\eps_{\cC}(c') - 
    \eps_{\cC}(c)g\eps_{\cC}(c'),
   $$
for all $c,c'\in \cC$,  and the unit is $g$. Furthermore 
 the counit $\eps_{\cC}:\cC\to R$ is a ring map split by ring maps 
 $i_{L},i_{R}: R\to \cC$ given by, $i_{L}: 
 r\mapsto rg$, $i_{R}:r\mapsto gr$. In the case of the canonical 
 Sweedler coring, the above ring structure 
 appears in 
 \cite[Section~1.2]{Nus:non} in the context of a braiding related to the
 noncommutative descent theory.

Another basic property of corings with a grouplike element is the
fact that their dual rings are augmentation rings. Recall from
\cite{Swe:pre} that the
left dual of $\cC$, ${}^*\cC = \lhom R \cC R$ 
is a unital ring with unit $\eps_{\cC}$ and the 
product 
    defined for all $\xi,\xi'\in {}^{*}\cC$ and $c\in\cC$ by
    $(\xi\xi')(c) = \xi(c\sw 1 \xi'(c\sw 2))$. 
    ${}^*\cC$ is an $(R,R)$-bimodule via $(r\xi r')(c) = \xi(cr)r'$. In
the case of the canonical coring associated to a ring extension $S\to
R$, ${}^*\cC$ is simply the left $S$-endomorphism ring of $R$, while in
the case of the $A$-coring $\cC = A\otimes C$ associated to an entwining
structure $(A,C)_\psi$, ${}^*\cC$ is a $\psi$-twisted convolution
product algebra  $({\rm Hom}(C,A), *_\psi)$  introduced in
\cite{Brz:coh} (see \cite{Kop:twi} for a special case).   
In general, every right
$\cC$-comodule $\cM$ is a left ${}^*\cC$-module via $\xi m = m\sw
0\xi(m\sw 1)$. This makes $\M^\cC_R$ a subcategory of the category of
left ${}^*\cC$-comodules. Recently it has been shown in \cite{Wis:com}
that $\M^\cC_R$ is a full subcategory of ${}_{^*\cC}\M$ if and only if
$\cC$ is a locally projective left $R$-module (cf.\ \cite{Zim:pur} for
the definition of a locally projective module).  Next recall
from \cite[p.\ 143]{CarEil:hom} that a ring $R$ is called a  {\em left 
augmentation ring} if there exists an $R$-module $M$ 
and a  left 
$R$-module morphism $\pi: R\to M$. $M$ is called an {\em augmentation 
module}. 

\begin{proposition}
    Let $\cC$ be an $R$-coring with a grouplike element $g$. Then 
${}^*\cC$  is a  left augmentation
 ring 
    with the augmentation module $R$. The left action of 
${}^*\cC$ on $R$ is provided by 
    $\xi r  = \xi(gr)$ for all $r\in R$, $\xi\in {}^*\cC$.
\Label{prop.augmentation}
\end{proposition}
    \begin{proof}
	 That the map given above defines a left action of 
${}^*\cC$ on $R$ follows from the fact that $R$ is a 
right $\cC$-comodule, and hence it is a left ${}^*\cC$-module as
recalled above. 
The augmentation 
	$\pi:{}^*\cC\to R$ is given by $\xi\mapsto \xi(g)$. 
	To show that $\pi$ is a 
	left ${}^*\cC$-module map, take any $\xi,\xi'\in {}^*\cC$ 
	and compute: $
	\pi(\xi\xi') = (\xi\xi')(g) = \xi(g \xi'(g)) = \xi(g\pi(\xi'))
	 = \xi\pi(\xi')$.
    \end{proof}

Proposition~\ref{prop.augmentation} extends an observation made in 
\cite{Kle:dua} in the case of finitely generated corings.

\section{The Amitsur complex}
\Label{section.differential}
There is a close relationship between corings with a grouplike 
element 
and a certain type of differential graded rings, first described by 
Rojter \cite{Roj:mat}.

First recall that
a {\em differential graded ring} is an $\N\cup\{ 0\}$-graded ring 
$\Omega = \bigoplus_{n=0}^{\infty}\Omega^{n}$ together with an 
additive degree 
one operation $d: \Omega^{\bullet}\to \Omega^{\bullet+1}$ such that 
$d\circ d =0$, which
satisfies the graded Leibniz rule, i.e.,  
for all elements $\omega'\in \Omega$ and all degree $n$ elements $\omega$,
$d(\omega\omega') = d(\omega)\omega' + (-1)^{n}\omega d(\omega')$.

Given a coring with a grouplike element one  constructs an associated 
differential
graded  ring as follows \cite{Roj:mat}.
 Consider a tensor ring
    $\Omega(\cC) = \bigoplus_{n=0}^{\infty}\Omega^{n}(\cC)$,
    where $\Omega^{0}(\cC) = R$ and $\Omega^{n}(\cC) = 
    \cC\otimes_{R}\cC\otimes_{R}\cdots\otimes_{R}\cC$ ($n$-times). 
    Define a degree one additive map $d:\Omega(\cC)\to\Omega(\cC)$ via
    $ d(r) = gr - rg$, and
    \begin{eqnarray*}
	d(c^{1}\otimes_{R}\cdots \otimes_{R}c^{n}) &=& 
	g\otimes_{R}c^{1}\otimes_{R}\cdots \otimes_{R}c^{n}\\
	&&\hspace{-1in}+ \sum_{i=1}^{n}(-1)^{i}c^{1}\otimes_{R}\cdots \otimes_{R}
	c^{i-1}\otimes_{R}\Delta_{\cC}(c^{i})\otimes_{R}c^{i+1} 
	\otimes_{R}\cdots\otimes_{R}c^{n} \\
	&&\hspace{-1in}+ (-1)^{n+1}c^{1}\otimes_{R}\cdots \otimes_{R}c^{n}\otimes_{R}g
    \end{eqnarray*}
A direct computation shows that $d\circ d =0$ and that $d$ satisfies the
Leibniz rule, hence $(\Omega(\cC),d)$ is a differential graded ring. 

\begin{remark}
In fact the above construction of a differential graded ring 
 can be extended to a slightly more
general situation. An element $g\in\cC$ is called a {\em
semi-grouplike element}\index{semi-grouplike element} provided $\DC(g) =
g\otimes_R g$. Note that every coring has a semi-grouplike element
(indeed, take $g=0$). Note also that the definition of a semi-grouplike
element implies that $u=\eC(g)\in R$ is an idempotent in the
centraliser of $g$ in $R$, i.e., $u^2 =u$ and $ug=gu$. The above
construction of a differential graded ring  makes no use 
of the fact that $u=1_R$ for
a grouplike element. Thus it holds even if
the phrase ``$g$ is a grouplike element'' is replaced by the phrase ``$g$ is
a semi-grouplike element''. The case of $g=0$ is of particular interest
for the development of the Cartier cohomology of corings \cite{Guz:coi}
as it corresponds to the cobar resolution of $\cC$.
\Label{remark.Amitsur}
\end{remark}
Following Sweedler \cite{Swe:gro}, given a 
ring $S$, a ring $R$ is called an {\em $S$-ring}\index{$S$-ring} or 
an {\em algebra over $S$}\index{algebra over a ring}
 if there is 
a ring map $S\to R$. Thus the notion of an $S$-ring  extends the
notion of an $S$-algebra to the case in which $S$ is neither commutative
nor central in $R$.
Extending further this notion to differential graded 
rings one says that $\Omega$ is an {\em $S$-relative differential 
graded ring}\index{relative differential graded ring}\index{$S$-relative differential graded
ring}\index{ring, relative differential graded} or a {\em 
differential graded algebra over $S$}\index{differential graded 
algebra over a ring}\index{algebra over a ring, graded differential} 
if there is  a ring map $S\to R=\Omega^{0}$ such that 
$d$ 
is an $(S,S)$-bimodule map and $d(S)=0$. 
\begin{proposition}
    Let $\cC$ be an $R$-coring with a grouplike element $g$, and let 
    $S=R^{co\cC}_{g}$. Then the differential graded ring 
$(\Omega(\cC),d)$ 
     is an $S$-relative 
    differential graded ring.
\end{proposition}
\begin{proof}
    This proposition immediately follows from the facts that $S$ is a 
    centraliser of $g$ in $R$ (hence $d(S)=0$) and that $\Delta_{\cC}$ 
    is an $(S,S)$-bimodule map.
\end{proof}
\begin{example}
    Take a ring extension $S\to R$, the canonical Sweedler's coring 
    $\cC = R\otimes_{S}R$ and a grouplike $g = 
1_{R}\otimes_{S}1_{R}$. 
    Then
    $
    \Omega^{n}(\cC) = R\otimes_{S}R\otimes_{S}\ldots\otimes_{S}R
    $
    ($n+1$-times), and $d^{n}= \sum_{i=0}^{n+1}(-1)^ie^{n}_{i}$, where
    $$
    e^{n}_{i}:r_{1}\otimes_{S}\ldots\otimes_{S} r_{n+1}\mapsto 
    r_{1}\otimes_{S}\ldots\otimes_{S}r_{i}\otimes_{S}1_{R}
    \otimes_{S}r_{i+1}\otimes_{S}\ldots\otimes_{S}r_{n+1},
    $$
    $i=0,1,\ldots, n+1$. This means that $(\Omega(R\otimes_{S}R),d)$ 
    is the Amitsur complex\index{Amitsur complex} associated to a ring 
    extension $S\to R$ \cite{Ami:sim} 
(see also \cite[Section~6]{Art:Azu}).
\Label{example.Amitsur} 
\end{example}

Motivated by Example~\ref{example.Amitsur} we call the cochain 
complex $(\Omega(\cC),d)$
the {\em Amitsur complex}
 associated to a coring $\cC$ and a 
grouplike $g$.
\begin{example}\Label{ex.psin}
    Let $\cC = A\otimes C$ be a coring associated to an entwining 
    structure
    $(A,C)_{\psi}$ over $k$. Suppose that $A$ is an 
    $(A,C)_{\psi}$-entwined module. Then
    $\Omega^{n}(\cC) = A\otimes C^{\otimes n}$, and
    \begin{eqnarray*}
    d^{n}(a\otimes c_{1}\otimes\ldots\otimes c_{n}) &=& a\sw 0\otimes 
    a\sw 1\otimes c_{1}\otimes\ldots\otimes c_{n}\\
    &&+ \sum_{i=1}^{n}(-1)^{i}a\otimes 
    c_{1}\otimes\ldots\otimes\Delta(c_{i}) \otimes \ldots\otimes 
    c_{n}\\
    &&+(-1)^{n+1}a\psi^n(c_{1}\otimes\ldots\otimes
c_{n}\otimes {1\sw 0})\otimes
1\sw 1,
    \end{eqnarray*}
where $\psi^n = (\psi\otimes C^{\otimes n-1})\circ \ldots 
\circ (C^{\otimes n-2}\otimes\psi\otimes
C)\circ (C^{\otimes n-1}\otimes \psi)$.
    Note that the product in $\Omega(\cC)$ reads
    \begin{eqnarray*}
   && (a\otimes c_{1}\otimes\ldots\otimes c_{m})(a'\otimes 
    c_{m+1}\otimes\ldots\otimes c_{m+n})\\
&& \hspace{3cm}= 
    a\psi^m(c_{1}\otimes\ldots\otimes c_{m}\otimes a')\otimes 
    c_{m+1}\otimes \ldots\otimes c_{m+n}.
    \end{eqnarray*}
\end{example}
\begin{proof}
Only the last term in the expression for $d^n$ might require some
explanation. We use the {\em $\alpha$-notation} $\psi(c\otimes a) =
a_\alpha\otimes c^\alpha$ (summation over repeated indices understood) so that
$\psi^n(c_1\otimes\ldots \otimes c_n\otimes a) = a_{\alpha_n\ldots
\alpha_1}\otimes c_1^{\alpha_1}\otimes\ldots \otimes c_n^{\alpha_n}$.
Then the last term in the expression for $d^n$ is obtained by 
the following chain of
identifications
\begin{eqnarray*}
&&(a\otimes c_1\otimes_A 1_A\otimes c_2\otimes_A\ldots \otimes_A
1_A\otimes c_n)\otimes_A 1\sw 0\otimes 1\sw 1 \\
&&~~~~~~~~~~~~~~~~~~~~=a\otimes c_1\otimes_A 1_A\otimes c_2\otimes_A\ldots \otimes_A
{1\sw 0}_{\alpha_n}\otimes c_n^{\alpha_n}\otimes 1\sw 1 \\
&&~~~~~~~~~~~~~~~~~~~~= a\otimes c_1\otimes_A 1_A\otimes c_2\otimes_A\ldots \otimes_A 
{1\sw 0}_{\alpha_n\alpha_{n-1}}\otimes c_{n-1}^{\alpha_{n-1}}
\otimes c_n^{\alpha_n}\otimes 1\sw 1\\
&&~~~~~~~~~~~~~~~~~~~~\ldots\\
&&~~~~~~~~~~~~~~~~~~~~= a{1\sw 0}_{\alpha_{n}\ldots\alpha_{1}}\otimes 
    c_{1}^{\alpha_{1}}\otimes\ldots\otimes c_{n}^{\alpha_{n}}\otimes
1\sw 1.
\end{eqnarray*}
Similar chain leads to the product in $\Omega(\cC)$.
\end{proof}

Let $\cC$ be an $R$-coring with a grouplike element $g$, and let
$S=R^{co\cC}_g$. Recall from \cite{Brz:str} that a pair 
$(\cC,g)$ is said to be a {\em Galois coring}
 if and only if
 there exists an $R$-coring
isomorphism $\chi: R\otimes_SR\to \cC$ such that 
$\chi(1_{R}\otimes_S1_{R})
=g$. The map $\chi$ is called a {\em Galois isomorphism}.
\begin{proposition}
    The Amitsur complex of a Galois $R$-coring $(\cC,g)$ is acyclic 
    provided $R$ is a faithfully flat left module of the subring $S$ of 
    its $g$-coinvariants.
\Label{prop.Amitsur.acyclic}
\end{proposition}
\begin{proof}
    Since $R$ is a faithfully flat left $S$-module and the Amitsur 
    complex $(\Omega(\cC),d)$ is a complex in the category of right 
    $S$-modules, it suffices to show that the complex 
    $(\Omega(\cC)\otimes_{S}R, d\otimes_{S}R)$ is acyclic. Let $\chi: 
    R\otimes_{S}R\to\cC$ be the  Galois isomorphism of 
    corings and write $\chi^{-1}(c) = c\su 1\otimes_{S}c\su 2$ 
    (summation understood), for all $c\in\cC$. First note that for 
    all $c\in\cC$, $c\su 1 gc\su 2 =c$, and then compute
    \begin{eqnarray*}
    &&(\cC\otimes_{R}\chi)(c\su 1g\otimes_{R}1_{R}\otimes_{S}c\su 2 - 
    c\sw 1\otimes_{R}c\sw 2\su 1\otimes_{S}c\sw 2\su 2) \\
    &&= c\su 
    1g\otimes_{R}gc\su 2-c\sw 1\otimes_{R}c\sw 2
    = \Delta_{\cC}(c\su 
    1gc\su 2 - c) = 0.
    \end{eqnarray*}
    Since $\chi$ is bijective we conclude that for all $c\in \cC$,
    $$
    c\su 1g\otimes_{S}c\su 2 = c\sw 1\chi^{-1}(c\sw 2).
    \eqno{(*)}
    $$
    
    For any $n=1,2,\ldots$ consider an additive map
    $h^{n}:\Omega^{n}(\cC)\otimes_{S}R\to 
    \Omega^{n-1}(\cC)\otimes_{S}R$, $
    h^{n}:c^{1}\otimes_{R}\ldots\otimes_{R}c^{n}\otimes_{S}r\mapsto 
    (-1)^{n}c^{1}\otimes_{R}\ldots\otimes_{R}c^{n-1}\chi^{-1}(c^{n}r)$.
    We will show that the collection $h$ of all such $h^{n}$ is a 
    contracting homotopy for $d\otimes_{S}R$.
    On  one hand we have
    \begin{eqnarray*}
	h^{n+1}(d^{n}(c^{1}\otimes_{R}\ldots\otimes_{R}c^{n})
	\otimes_{S}r) &=& 
	(-1)^{n+1}g\otimes_{R}c^{1}\otimes_{R}\ldots\otimes_{R}c^{n-1} 
	\chi^{-1}(c^{n}r)\\
	&& \hspace{-1.8in}+\sum_{i=1}^{n-1}(-1)^{n+i+1}c^{1}\otimes_{R}\ldots \otimes_{R}
	\Delta_{\cC}(c^{i})\otimes_{R}\ldots\otimes_{R}c^{n-1}
	\chi^{-1}(c^{n}r) \\
	&&\hspace{-1.8in}- c^{1}\otimes_{R}\cdots \otimes_{R}c^{n}\sw 
	1\chi^{-1}(c^{n}\sw 2r)
	+ c^{1}\otimes_{R}\cdots \otimes_{R}c^{n}\chi^{-1}(gr)\\
	&&\hspace{-1.9in}=
	(-1)^{n+1}g\otimes_{R}c^{1}\otimes_{R}\ldots\otimes_{R}c^{n-1} 
	\chi^{-1}(c^{n}r)\\
	&&\hspace{-1.8in} +\sum_{i=1}^{n-1}(-1)^{n+i+1}c^{1}\otimes_{R}\ldots \otimes_{R}
	\Delta_{\cC}(c^{i})\otimes_{R}\ldots\otimes_{R}c^{n-1}
	\chi^{-1}(c^{n}r) \\
	&&\hspace{-1.8in}- c^{1}\otimes_{R}\cdots \otimes_{R}{c^{n}}\su
	1g\otimes_{S}{c^{n}}\su 2r
	+ c^{1}\otimes_{R}\cdots \otimes_{R}c^{n}\otimes_{S}r,
    \end{eqnarray*}
    where we have used equation ($*$) and the fact that 
    $\chi^{-1}(g) = 1_{R}\otimes_{S}1_{R}$. On the other hand,
    \begin{eqnarray*}
	d^{n-1}(h^{n}(c^{1}\otimes_{R}\ldots\otimes_{R}c^{n})
	\otimes_{S}r) &=& 
	(-1)^{n}g\otimes_{R}c^{1}\otimes_{R}\ldots\otimes_{R}c^{n-1} 
	\chi^{-1}(c^{n}r)\\
	&& \hspace{-1.5in}+\sum_{i=1}^{n-1}(-1)^{n+i}c^{1}\otimes_{R}\ldots \otimes_{R}
	\Delta_{\cC}(c^{i})\otimes_{R}\ldots\otimes_{R}c^{n-1}
	\chi^{-1}(c^{n}r) \\
	&&\hspace{-1.5in}+ c^{1}\otimes_{R}\cdots \otimes_{R}{c^{n}}\su
	1g\otimes_{S}{c^{n}}\su 2r.
    \end{eqnarray*}
    Thus $h^{n+1}d^{n}+d^{n-1}h^{n}= \Omega^{n}(\cC)\otimes_{S}R$. 
    This means that $d\otimes_{S}R$ is homotopic to the identity, so 
    the complex $(\Omega(\cC)\otimes_{S}R, d\otimes_{S}R)$ is 
    acyclic, and therefore the Amitsur complex is acyclic by the 
    virtue of the fact that the functor $-\otimes_{S}R$ reflects 
    exact sequences (for $R$ is a faithfully flat left $S$-module).
    \end{proof}

    In fact, since there is a ring map $S\to R$, hence a map $S\to 
    \Omega(\cC)$, if $S\to R$ is faithfully flat and $(\cC,g)$ is a 
    Galois coring then the associated Amitsur 
    complex $\Omega(\cC)$ is a resolution of $S$. Therefore 
    Proposition~\ref{prop.Amitsur.acyclic} can be seen as a coring 
    version of a well-known fact in the faithfully flat descent 
    theory (cf.\ \cite[Proposition~6.2]{Art:Azu}).

\section{Connections and comodules}
\Label{sec.connection.grouplike}
The category of comodules of a coring with a grouplike element can be 
naturally described in terms of connections over a ring. 
   Let $S\to R$ be a ring extension, and let $\Omega$ be an 
$S$-relative 
    differential graded ring with $R=\Omega^{0}$. Recall from
\cite{Con:dif} that a {\em
connection} 
    in a right $R$-module $M$ is a right $S$-linear map $\nabla: 
    M\otimes_{R}\Omega^{\bullet}\to M\otimes_{R}\Omega^{\bullet+1}$,
which satisfies the Leibniz rule, i.e., 
    such that for all $\omega\in M\otimes_{R}\Omega^{k}$ and 
    $\omega'\in\Omega$,
    $\nabla(\omega\omega') = \nabla(\omega)\omega'+(-1)^{k}\omega 
    d(\omega')$.
    A {\em curvature}
    of a connection $\nabla$ is a right $S$-linear map 
    $F_{\nabla}:M\to M\otimes_{R}\Omega^{2}$,
    defined as a restriction of $\nabla\circ\nabla$ to $M$, i.e., 
    $F_{\nabla}= \nabla\circ\nabla\mid_{M}$. A connection is said to 
    be {\em flat}
 if its curvature is identically equal to 0.

It is important to note that a connection is fully determined by its 
restriction to the module $M$. Indeed, any element of 
$M\otimes_{R}\Omega$ is a sum of terms of the form 
$m\otimes_{R}\omega$ with $m\in M$ and $\omega\in \Omega$. Now, 
viewing $m$ as an element of $M\otimes_{R}\Omega^{0}$, and using 
the Leibniz rule, the action of $\nabla$ on $m\otimes_{R}\omega$ 
reads,
$\nabla(m\otimes_{R}\omega) = \nabla(m)\omega +m\otimes_{R}d(\omega)$.

 To describe a relationship between connections and comodules of a 
coring we first need to introduce an appropriate differential graded 
ring \cite[Lemma~1]{Roj:mat}.
\begin{proposition}
    Let $\cC$ be an $R$-coring with a grouplike $g\in \cC$ and let 
$S$ 
    be the subring of $g$-coinvariants. Then  the associated Amitsur 
    complex $(\Omega(\cC),d)$ restricts to the $S$-relative 
    differential graded ring $(\Omega(\cC/S),d)$ with $\Omega^{0}(\cC/S) 
= 
    R$ and
    $$
    \Omega^{n}(\cC/S) = 
    \ker\eps_{\cC}\otimes_{R}\ker\eps_{\cC}\otimes_{R}
    \ldots\otimes_{R}\ker\eps_{\cC},
    $$
    ($\ker\eps_{\cC}$ taken $n$-times).
\end{proposition}
\begin{proof}
The key observation here is that, firstly, for all $r\in R$, 
$\eps_{\cC}(d(r)) = \eps_{\cC}(gr)-\eps_{\cC}(rg) = r-r =0$, and, 
secondly, for any $c^{1},\ldots, c^{n}\in\ker\eps_{\cC}$, the Amitsur 
coboundary operator $d^{n}$ can be written equivalently as
\begin{eqnarray*}
    d^{n}(c^{1}\otimes_{R}\ldots \otimes_{R}c^{n}) &=& 
    \sum_{i=1}^{n}(-1)^{i}c^{1}\otimes_{R}\ldots\otimes_{R}c^{i-1}\\
    &&\otimes_{R}(c^{i}\sw 1-g\eps_{\cC}(c^{i}\sw 
    1))\otimes_{R}(c^{i}\sw 2 - \eps_{\cC}(c^{i}\sw 
    2)g)\\
    &&\otimes_{R}c^{i+1}\otimes_{R}\ldots\otimes_{R}c^{n}.
\end{eqnarray*}
This expression shows immediately that the image of 
$d^{n}$ applied to $(\ker\eps_{\cC})^{\otimes_{R}n}$ is  
in $(\ker\eps_{\cC})^{\otimes_{R}n+1}$ as required.
\end{proof}

Thus we have constructed an $S$-relative differential graded ring 
$\Omega(\cC/S)$, which is termed a ring of {\em $\cC$-valued 
differential forms on $R$}. This construction has an interesting 
converse. Let $\Omega$ be a 
 differential graded ring   with $\Omega^{0}=R$, and $\Omega^{n}= 
 \Omega^{1}\otimes_{R}\Omega^{1}\otimes_{R}\ldots\otimes_{R}\Omega^{1}$ 
 (n-times). Consider  a left $R$-module  
 $\cC= Rg \oplus \Omega^{1}$. One easily checks that $\cC$ can be made 
 into an $(R,R)$-bimodule with a right multiplication
$ (rg+\omega)r' = rr'g+ rd(r') +\omega r'$ for all $r,r'\in R$, 
$\omega\in \Omega^{1}$. 
Then $\cC$ is an $R$-coring with the coproduct and counit
$$
\DC (rg) = rg\otimes_{R} g, \quad \DC(\omega) = g\otimes_{R}\omega 
+\omega\otimes_{R}g - d(\omega), \quad \eC(rg+\omega) = r,
$$
for all $r\in R$ and $\omega \in \Omega^{1}$. The coassociativity
of $\DC$ follows from the equality $d(\omega\su 1)\otimes_R \omega\su 2 +
\omega\su 1\otimes_R d(\omega\su 2) = 0$, where 
$d(\omega) = \omega\su 1\otimes_{R}\omega\su 2$ is a notation. This 
is a consequence of the
Leibniz rule and the nilpotency of $d$. Clearly, $g$ 
is a grouplike element and $\Omega = \Omega(\cC/S)$, where 
 $S = \ker(d: R\to \Omega^{1})$. Note also that $d(r) = gr-rg$. 
 In particular some differential
calculi of noncommutative geometry based on Dirac operators or Fredholm
modules, lead to corings with grouplike elements.

Given a ring extension $S\to R$ there is an associated universal 
differential
graded algebra over $S$ generated by the $S$-ring $R$ known as
 an
{\em algebra of $S$-relative differential forms} 
$\Omega_SR$\index{relative
differential forms} \cite{CunQui:alg}. Let $R/S$ denote the cokernel of
the ring homomorphism $S\to R$ viewed as a map of $(S,S)$-bimodules, and
let $\pi: R\to R/S$ be the canonical map of $(S,S)$-bimodules.
Thus $R/S$ is an $S$-bimodule and for any $n\in \N\cup\{0\}$ one can 
consider an $(R,S)$-bimodule 
$$
\Omega^n_S R = R\otimes_S(R/S)^{\otimes_Sn} = R\otimes_S R/S\otimes_S
R/S\otimes_S\cdots \otimes_S R/S,
$$
and combine them into a direct sum $\Omega_SR = \oplus_{n=0}^\infty
\Omega^n_S R$. $\Omega_SR$ is a cochain complex with a coboundary
operator 
$$
d: \Omega^\bullet_SR\to \Omega^{\bullet+1}_SR, \qquad r_0\otimes_S
r_1\otimes_S\ldots \otimes_S r_n\mapsto 1_R\otimes_S\pi(r_0)\otimes_S
r_1\otimes_S \ldots \otimes_Sr_n.
$$
It is clear that $d\circ d = 0$ since $S\to R$ as a ring homomorphism
is a unit preserving map so that $\pi(1_R) = 0$. Less trivial is the
observation that $(\Omega_SR, d)$ is an $S$-relative differential
graded ring. The product in $\Omega_SR$ is given by the formula
$$
(r_0,\ldots, r_n)(r_{n+1},\ldots, r_m) = \sum_{i=0}^{n}(-1)^{n-i}
(r_0,\ldots, r_{i-1}, r_i\cdot r_{i+1}, r_{i+2},\ldots ,r_m),
$$
where we write $(r_0,\ldots, r_n)$ for $r_0\otimes_S\ldots\otimes_S r_n$
etc., to relieve the notation.  Also, the notation $r_i\cdot r_{i+1}$ is
a formal expression that is to be understood as follows. Take any 
$r'_i\in \pi^{-1}(r_i)$ and $r'_{i+1} \in \pi^{-1}(r_{i+1})$,
$i=1,\ldots n-1$ then
$$
r_i\cdot r_{i+1} = \left\{ \begin{array}{ll}
			r_ir'_{i+1} 	   & \mbox{for $i=0$}\\
			\pi(r'_i r'_{i+1}) &\mbox{for $0<i<n$}\\
			\pi(r'_ir_{i+1})   & \mbox{for $i=n$}  
			\end{array}
			\right.
$$
Although $r_i\cdot r_{i+1}$ depends on the choice of the $r'_i$, it can
be easily shown that the product of cochains does not. It can also be
shown that the above expression defines an associative product and that
$d$ satisfies the graded Leibniz rule. The details can be found in
\cite{CunQui:alg}. The algebra of $S$-relative differential forms has
the following universality property. Given any graded differential ring
$\Omega = \oplus_n \Omega^n$ and a ring homomorphism $u: R\to \Omega^0$
such that $d(u(s))= 0$ for all $s\in S$, there exists a unique
differential graded ring homomorphism $u_* :\Omega_SR\to \Omega$
extending $u$. This means, in particular, that the identity map $R\to R$
extends to a map of $S$-relative differential graded rings $\Omega_SR\to
\Omega(\cC/S)$. This can be understood purely in terms of coring valued 
differential forms since we have the following 
\begin{proposition}
Let $S\to R$ be  a ring extension. Take the canonical Sweedler's coring
$\cC = R\otimes_S R$, and a grouplike element $g =1_R\otimes_S1_R$. Then
the differential graded algebra over $S$ of $\cC$-valued differential 
forms is isomorphic to the algebra of $S$-relative differential forms 
$\Omega_SR$.
\Label{prop.cuntz} 
\end{proposition}
\begin{proof}
We first explicitly describe the structure of $\cC = R\otimes_SR$-valued
differential forms. Note that since the counit of the canonical coring
coincides with the product map $\mu_{R/S} :R\otimes_S R\to R$, $r\otimes_S
r'\mapsto rr'$, we have $\Omega^1(\cC/S) = \ker\eps_\cC  = \ker\mu_{R/S}$.
Clearly, the Amitsur $0$-differential $d: r\mapsto 1_R\otimes_S r - r\otimes_S
1_{R}$, has values   restricted to $\ker\mu_{R/S}$. Note that 
$\Omega^n(\cC/S) = (\ker\mu_{R/S})^{\otimes_R n}$. Thus if we could show that
$\ker\mu_{R/S} \cong R\otimes_S R/S$ as $(R,R)$-bimodules,
 then we would have the required form of
$\Omega^n(\cC/S)$. Indeed, the following iteration
\begin{eqnarray*}
\Omega^n(\cC/S) &=& \Omega^{n-1}(\cC/S)\otimes_R\ker\mu_{R/S}\\
& \cong &
\Omega^{n-1}(\cC/S) \otimes_R R\otimes_S R/S\cong  
\Omega^{n-1}(\cC/S)\otimes_S R/S,
\end{eqnarray*}
repeated $n$-times would yield the desired result. 
Note that $\ker\mu_{R/S} \cong R\otimes_S R/S$
as $(R,S)$-bimodules via the map $\theta: \ker\mu_{R/S} \to R\otimes_S R/S$,
$\sum_i r_i\otimes_S r'_i\mapsto \sum_i r_i\otimes_S \pi(r'_i)$ with the
inverse $\theta^{-1}: r\otimes_S \pi(r') = r\otimes_S r' - rr'\otimes_S
1_R$. Note also that the map $\theta^{-1}$ does not depend on the choice of
$r'$ in the inverse image of $\pi(r')$. Indeed, if $\pi(r') = 0$ then
$r' = s1_R$ with $s\in S$ and hence $r\otimes_S r' - rr'\otimes_S 1_R= 
r\otimes_S s1_R - rs\otimes_S 1_R =0$ as needed. The right $R$-module
structure of $R\otimes_S
R/S$ is derived from the product in $\Omega_SR$, i.e., $(r_0\otimes_{S}
\pi(r_1)) r = -r_0r\otimes_S \pi(r_1) + r_0\otimes_S\pi(r_1r)$ and is
well defined (does
not depend on the choice of $r_1$) by the similar argument as above.
Clearly $\theta$ and $\theta^{-1}$ are maps of right $R$-modules.
The isomorphism $\theta$ extends to cochains of all degrees and 
one can easily check that it
provides an isomorphism of $S$-relative graded differential algebras. 
This isomorphism involves projections $\pi$ in all bar the first tensorand
and thus maps all the Amitsur operators $e_i^n$ with $i>0$ in
Example~\ref{example.Amitsur}  to  $0$. Thus the resulting differential 
$d$ has the
form $d:r_0\otimes_S
r_1\otimes_S\ldots \otimes_S r_n\mapsto 1_R\otimes_S\pi(r_0)\otimes_S
r_1\otimes_S\ldots \otimes_S r_n
$, as required, and the proof of the proposition is completed.
\end{proof}

Next we address the problem of existence of $\Omega(\cC/S)$-valued connections.

\begin{theorem}
Let $\cC$ be an $R$-coring with a grouplike element $g\in \cC$, let $S$ 
    be the subring of $g$-coinvariants of $R$, and let $\Omega(\cC/S)$ be the 
    ring of $\cC$-valued differential forms on $R$. 
    A right $R$-module admits an $\Omega(\cC/S)$ valued connection if 
and only if $M\otimes_{R}\eC$ is a right $R$-module retraction.
\Label{thm.ex.con}
\end{theorem}
\begin{proof}
    Given a connection $\nabla: M\to M\otimes_{R}\Omega^{1}(\cC/S)$ 
    define an additive map
    $$
    j_{\nabla}: M\to M\otimes_{R}\cC, \qquad m\mapsto \nabla(m) 
    +m\otimes_{R} g.
    $$
    Since ${\rm Im}(\nabla) \subseteq M\otimes_R\ker\eC$ and $\eC(g) 
    =1_{R}$, the map $j_{\nabla}$ is an additive section of 
    $M\otimes_{R} \eC$. Furthermore for all $m\in M$ and $r\in 
    R$,
    \begin{eqnarray*}
	j_{\nabla}(mr) &=& \nabla(mr) +mr\otimes_{R}g = \nabla(m)r 
    +m\otimes_{R}d(r) +mr\otimes_{R}g\\
    &=& \nabla(m)r +m\otimes_{R}gr - m\otimes_{R}rg +m\otimes_{R}rg
    = j_{\nabla}(m)r.
    \end{eqnarray*}
    where we used the Leibniz rule for $\nabla$ to obtain the second 
    equality. Thus $j_{\nabla}$ is a right $R$-linear section of 
    $M\otimes_{R} \eC$.
    
    Conversely, suppose $j: M\to M\otimes_{R} \cC$ is a right $R$-linear 
    section of $M\otimes_{R} \eC$, and define an additive map
    $$
    \nabla_{j}:M\to M\otimes_{R}\Omega^{1}(\cC/S) = 
    M\otimes_{R}\ker\eC, \qquad m\mapsto j(m) - m\otimes_{R} g.
    $$
    Note that $\nabla_{j}$ is well-defined since the fact that $j$ is 
    a section implies that for all $m\in M$, 
    $\sum_{i}m^{i}\eC(c^{i}) = m$, where 
    $\sum_{i}m^{i}\otimes_{R} c^{i} = j(m)$. Therefore
    $$
    \nabla_{j}(m) = 
    \sum_{i}(m^{i}\otimes_{R} c^{i}-m^{i}\eC(c^{i})\otimes_{R} g) = 
    \sum_{i}m^{i}\otimes_{R} (c^{i}-\eC(c^{i})g),
    $$
    and each $c^{i}-\eC(c^{i})g\in\ker\eC$. Finally  for all $m\in M$ and $r\in R$ we have
    \begin{eqnarray*}
	\nabla_{j}(mr) &=& j(mr) -mr\otimes_{R} g = j(m)r - m\otimes_{R}rg\\
	&=& j(m)r -m\otimes_{R} gr +m\otimes_{R} gr -m\otimes_{R}rg
	= \nabla_{j}(m)r +m\otimes_{R} dr,
    \end{eqnarray*}
     so that $\nabla_{j}$ satisfies the Leibniz rule for a 
    connection as required.
\end{proof}

In particular Theorem~\ref{thm.ex.con} implies that if a right 
$R$-module $M$ admits an $\Omega(\cC/S)$ valued connection then $M$ is 
a direct summand of a right $R$-module $M\otimes_{R}\cC$. As a special 
case of Theorem~\ref{thm.ex.con} one can take an algebra $A$ over 
a field $k$, view it as an extension of $k$, and thus consider the 
corresponding Sweedler $A$-coring $\cC = A\otimes A$ and a grouplike 
$g=1_{A}\otimes 1_{A}$. In this case the $g$-coinvariants of $A$ are 
simply equal to $k$, and $\Omega(\cC/k) = \Omega A$ is a universal 
differential envelope of $A$ by Proposition~\ref{prop.cuntz}. By
the canonical identification $M\otimes_{A}A\otimes A\cong M\otimes A$, 
for any $M\in \M_{A}$, the map $M\otimes_{A}\eC$ coincides with the 
action of $A$ on $M$, $m\otimes a\mapsto ma$. The action is an 
$A$-module retraction if and only if $M$ is a projective right 
$A$-module. Thus, as a special case of Theorem~\ref{thm.ex.con}, we 
obtain the classic result of Cuntz and Quillen 
\cite[Corollary~8.2]{CunQui:alg} which relates the existence of 
(universal) connections to projectivity.

The importance of coring valued differential forms in relation to
comodules of a coring is revealed by the following theorem.
\begin{theorem}
    Let $\cC$ be an $R$-coring with a grouplike $g\in \cC$, let $S$ 
    be the subring of $g$-coinvariants, and let $\Omega(\cC/S)$ be the 
    ring of $\cC$-valued differential forms on $R$. Then a right 
    $R$-module $\cM$ is a right $\cC$-comodule if and only if it admits 
    a flat connection $\nabla: \cM\to \cM\otimes_{R}\Omega(\cC/S)$.
\Label{theorem.flat}
\end{theorem}
\begin{proof}
    Suppose $\cM$ is a right $\cC$-comodule with a coaction $\rho^{\cM}$ 
    and define
    $$
    \nabla: \cM\to \cM\otimes_{R}\Omega^{1}(\cC/S), \qquad m\mapsto 
    \rho^{\cM}(m) - m\otimes_{R}g.
    $$
    Since $\rho^{\cM}$ is a right $R$-module splitting of 
    $\cM\otimes_{R}\eC$, the map $\nabla$ is a connection by 
    Theorem~\ref{thm.ex.con}.
    
    We can now compute the curvature $F_{\nabla}$ of $\nabla$. Take 
    any $m\in \cM$ and compute
    \begin{eqnarray*}
	F_{\nabla}(m) &=& \nabla(m\sw 0\otimes_{R}m\sw 1- m\otimes_{R}g)\\
	&=& \nabla(m\sw 0)\otimes_{R}m\sw 1+m\sw 0\otimes_{R}dm\sw 
	1-\nabla(m)\otimes_{R}g +m\otimes_{R}dg\\
	&=& 0,
    \end{eqnarray*}
    by the coassociativity of $\rho^{\cM}$ and the definition of $d$.
    
    Conversely, suppose $\cM$ is a right $R$-module with a flat 
    connection $\nabla:\cM\to \cM\otimes_{R}\ker\eps_{\cC}$. For any 
$m\in 
    \cM$, write $\nabla(m) = \sum_{i}m^{i}\otimes_{R}c^{i}$. Then the 
    flatness of $\nabla$ means that $
    0 = \nabla(\sum_{i}m^{i}\otimes_{R}c^{i}) = \sum_{i}\nabla(m^{i})
    \otimes_{R}c^{i}+ \sum_{i}m^{i}\otimes_{R}d(c^{i})$,
    i.e.,
    $$
    \sum_{i}m^{i}\otimes_{R}
    c^{i}\sw 1 \otimes_{R}c^{i}\sw 2 = 
\sum_{ij}m^{ij}\otimes_{R}\tilde{c}^{j}\otimes_{R}c^{i}+
    \sum_{i}m^{i}\otimes_{R}g\otimes_{R}c^{i}
    +\sum_{i}m^{i}\otimes_{R}c^{i}\otimes_{R}g,
    $$
    where
    $\nabla(m^{i}) = \sum_{j}m^{ij}\otimes_{R}\tilde{c}^{j}$.
    
    Now define an additive map $
    \rho^{\cM}:\cM\to \cM\otimes_{R}\cC$,  $m \mapsto \nabla(m) 
    +m\otimes_{R}g$.
    Note that $\rho^{\cM}$ coincides with the map $j_\nabla$
constructed in the proof of Theorem~\ref{thm.ex.con}, and thus 
it is a right $R$-module map. 
    Since every $c^{i}$ in $\nabla(m) = 
    \sum_{i}m^{i}\otimes_{R}c^{i}$ is an element of $\ker\eps_{\cC}$,
and $\eps_\cC(g) =1_R$,  
    we immediately  conclude that 
    $(\cM\otimes_{R}\eps_{\cC})\rho^{\cM}(m) = m$. Therefore it only 
remains to be
    shown that $\rho^{\cM}$ is coassociative. Since explicitly 
    $\rho^{\cM}(m) = \sum_{i}m^{i}\otimes_{R}c^{i} +m\otimes_{R}g$, we 
    have
    \begin{eqnarray*}
	(\rho^{\cM}\otimes_{R}\cC)\rho^{\cM}(m) &=& 
	\sum_{i}\rho^{\cM}(m^{i})\otimes_{R}c^{i}+ \rho^{\cM}(m)\otimes_{R}g\\
	&=& \sum_{ij}m^{ij}\otimes_{R}\tilde{c}^{j}\otimes_{R}c^{i}+
    \sum_{i}m^{i}\otimes_{R}g\otimes_{R}c^{i}\\
    &&+\sum_{i}m^{i}\otimes_{R}c^{i}\otimes_{R}g 
    +m\otimes_{R}g\otimes_{R}g\\
    &=& \sum_{i}m^{i}\otimes_{R}
    c^{i}\sw 1 \otimes_{R}c^{i}\sw 2+m\otimes_{R}g\otimes_{R}g
    =(\cM\otimes_{R}\Delta_{\cC})\rho^{\cM}(m),
    \end{eqnarray*}
    where we used the flatness of $\nabla$ to derive the penultimate
equality. This proves that $\rho^{\cM}$ is a coassociative coaction, and 
    hence $\cM$ is a right $\cC$-comodule.
\end{proof}

In the proof of Theorem~\ref{theorem.flat} we have constructed two 
assignments. Given a right $R$-module $\cM$, to every right 
$\cC$-coaction $\rho^{\cM}$ one assigns a flat connection 
$\nabla_{\rho^{\cM}}$ defined by $\nabla_{\rho^{\cM}}: m\mapsto 
\rho^{\cM}(m)- 
m\otimes_{R}g$. Conversely to any flat connection $\nabla$ one 
assigns 
a right $\cC$-coaction $\rho^{\cM}_{\nabla}: m\mapsto \nabla(m) + 
m\otimes_{R}g$. Clearly these assignments are inverses of each other 
and hence establish an isomorphism of sets  of flat 
connections and right $\cC$-comodule structures. 
 In fact Theorem~\ref{theorem.flat} describes 
 an isomorphism of categories.

Consider a ring extension $S\to R$ and an $S$-relative differential 
graded ring $\Omega$ with $\Omega^{0}=R$. A category of (right) 
connections with values in $\Omega$, denoted by $\con(R/S,\Omega)$ 
consists of pairs $(M,\nabla)$ where $M$ is a right $R$-module and 
$\nabla: M\otimes_{R}\Omega^{\bullet}\to 
M\otimes_{R}\Omega^{\bullet+1}$ is a connection. A morphism 
$(M,\nabla) \to (N,\nabla')$ in 
$\con(R/S,\Omega)$ is a right $R$-module map $f:M\to N$ making the 
following diagram
\begin{diagram}
    M & \rTo^{f} & N\\
    \dTo^{\nabla} & & \dTo_{\nabla'}\\
    M\otimes_{R}\Omega^{1} & \rTo^{f\otimes_{R}\Omega^{1}} & 
    N\otimes_{R}\Omega^{1}
\end{diagram}
commute. Note that the Leibniz rule renders  also the following 
diagram commutative
\begin{diagram}
    M \otimes_{R}\Omega^{\bullet}& \rTo^{f\otimes_{R}\Omega^{\bullet}} & N\otimes_{R}\Omega^{\bullet}\\
    \dTo^{\nabla} & & \dTo_{\nabla'}\\
    M\otimes_{R}\Omega^{\bullet+1} & 
    \rTo^{f\otimes_{R}\Omega^{\bullet+1}} & 
    N\otimes_{R}\Omega^{\bullet+1}
\end{diagram}
 In particular we can consider the following diagram
\begin{diagram}
    M & \rTo^{\nabla} & M\otimes_{R}\Omega^{1} & \rTo^{\nabla} & 
    M\otimes_{R}\Omega^{2}\\
    \dTo_{f} & & \dTo_{f\otimes_{R}\Omega^{1}} && 
    \dTo_{f\otimes_{R}\Omega^{2}}\\
    N & \rTo^{\nabla'} & N\otimes_{R}\Omega^{1} & \rTo^{\nabla'} & 
    N\otimes_{R}\Omega^{2}
\end{diagram}
in which both left and right squares commute. This implies that the 
outer rectangle is commutative, hence $F'\circ f = 
(f\otimes_{R}\Omega^{2})\circ F$, where $F$ is the curvature of 
$\nabla$ and $F'$ is the curvature of $\nabla'$. This shows that 
$\con(R/S,\Omega)$ contains a full subcategory $\con_{0}(R/S,\Omega)$ 
of flat connections. Objects of $\con_{0}(R/S,\Omega)$ 
 are pairs $(M,\nabla)$ where $M$ 
is a right $R$-module and $\nabla$ is a flat connection.

Now Theorem~\ref{theorem.flat} leads to the following

\begin{theorem}
    Let $\cC$ be an $R$-coring with a grouplike element $g\in \cC$, let $S$ 
    be the subring of $g$-coinvariants, and let $\Omega(\cC/S)$ be the 
    ring of $\cC$-valued differential forms on $R$. Then 
    $\M_{R}^{\cC}$ is isomorphic to $\con_{0}(R/S, \Omega(\cC/S))$.
\Label{thm.iso.flat}
\end{theorem}
\begin{proof}
    On objects, the isomorphism is provided by the assignment 
    constructed in Theorem~\ref{theorem.flat}, while on morphisms 
    $\ff :\cM\to \cN$, $\ff\leftrightarrow\ff$. Indeed, if $\ff$ is a 
    morphism of right $\cC$-comodules then for all $m\in\cM$,
    \begin{eqnarray*}
	\nabla_{\rho^{\cN}}\circ \ff(m) &=& \rho^{\cN}\circ \ff(m) - 
	\ff(m)\otimes_{R}g \\
	&=& (\ff\otimes_{R}\cC)(\rho^{\cM}(m) - m\otimes_{R}g) = 
	(\ff\otimes_{R}\cC)\circ \nabla_{\rho^{\cM}}(m).
    \end{eqnarray*}
    Conversely, if $\ff$ is a morphism $(\cM,\nabla)\to (\cN,\nabla')$ 
    in $\con_{0}(R/S,\Omega(\cC/S))$ then
    \begin{eqnarray*}
	\rho^{\cN}_{\nabla'}\circ\ff(m) &=& \nabla'(\ff (m)) + 
	\ff(m)\otimes_{R} g = (\ff\otimes_{R}\cC)\circ\nabla(m) + 
	\ff(m)\otimes_{R}g\\
	&=& (\ff\otimes_{R}\cC)\circ(\nabla(m) + 
	m\otimes_{R}g) = (\ff\otimes_{R}\cC)\circ\rho^{\cM}_{\nabla}.
    \end{eqnarray*}
    This completes the proof of the theorem.
\end{proof}

\begin{example}
    The ring $R$ admits a flat connection, $\nabla(r) = gr-rg=d(r)$.
\end{example}
\begin{example}
    Let $(A,C)_{\psi}$ be an entwining structure over $k$, and let 
    $\cC=A\otimes C$ be the associated $A$-coring. Suppose $A$ is an 
    $(A,C)_{\psi}$-entwined module. Then $A\otimes C^{\otimes n}$, $n 
    =1,2,\dots$ has a flat connection
    $$
    \nabla(a\otimes c_{1}\otimes\ldots\otimes c_{n}) = a\otimes 
    c_{1}\otimes\ldots\otimes \Delta(c_{n}) - a\psi^n(c_{1}
    \otimes\ldots\otimes c_{n} \otimes {1\sw 
    0})\otimes 1\sw 1,
    $$
where $\psi^n : C^{\otimes n}\otimes A\to A\otimes C^{\otimes n}$ is
defined in Example~\ref{ex.psin}.
\Label{ex.con.ac}
\end{example}
\begin{proof}
This follows immediately from the fact that $A\otimes C^{\otimes n}$ is
a right $\cC$-comodule (or, equivalently, an $(A,C)_\psi$-entwined
module) \cite{Brz:mod}.
\end{proof}

Dually we have
\begin{example}
    Let $(A,C)_{\psi}$ be an entwining structure over $k$, let 
    $\cC=A\otimes C$ be the associated $A$-coring. Suppose $A$ is an 
    $(A,C)_{\psi}$-entwined module. Then $C\otimes A^{\otimes n}$, $n 
    =1,2,\dots$ has a flat connection
    $$
    \nabla(c\otimes a^{1}\otimes\ldots\otimes a^{n}) = c\sw 1\otimes 
    \psi_n(c\sw 2 \otimes a^{1}\otimes\ldots\otimes a^{n}) - c\otimes a^{1}\otimes\ldots 
    \otimes a^{n}1\sw 0\otimes 1\sw 1,
    $$
where $\psi_n= (A^{\tens n-1}\otimes \psi)\circ \ldots \circ (A\otimes
\psi\otimes A^{\tens n-2})\circ  (\psi\otimes A^{\tens n-1}): 
C\otimes A^{\tens n}\to A^{\tens n}\otimes C$.
\Label{ex.con.ca}
\end{example}

\begin{center}
{\sc Acknowledgements}
\end{center}
I would like to thank  
EPSRC for an Advanced Research Fellowship.

\end{document}